\documentstyle[12pt, fleqn]{article}
\renewcommand{\baselinestretch}{1.3}
\newcommand{\be}{\begin{equation}}
\newcommand{\ee}{\end{equation}}
\newcommand{\ba}{\begin{array}}
\newcommand{\ea}{\end{array}}
\newcommand{\bea}{\begin{eqnarray}}
\newcommand{\eea}{\end{eqnarray}}
\newcommand{\su}{\mbox{$SU_q(2)$\,\,}}
\newcommand{\pro}{\partial}
\newcommand{\op}{\omega^{++}}
\newcommand{\om}{\omega^{--}}
\newcommand{\oo}{\omega^0}
\newcommand{\dfrac}{\displaystyle\frac}
\newcommand{\dpst}{\displaystyle}
\newcommand{\dpp}{{\cal D}^{++}}
\newcommand{\dmm}{{\cal D}^{--}}
\newcommand{\docal}{{\cal D}^{0}}
\newcommand{\got}{\cal}
\newcommand{\vareps}{\varepsilon}
\newcommand{\rmm}{({\hat R}^{-1})}
\newcommand{\qin}{\displaystyle\frac{1}{q}}
\newcommand{\uone}{\mbox{$U(1)$\,\,}}
\newcommand{\uqtwo}{\mbox{$U_q(2)$\,\,}}
\newcommand{\suu}{\mbox{$SU(2)$\,}}
\newcommand{\nn}{\nonumber}

\begin{document}
\title{
             REPRESENTATION OF $SU_q(2)$-COVARIANT $q$-LIE ALGEBRA
              IN TERMS OF DIFFERENTIAL OPERATORS }

\author{ D.G. Pak \thanks {E-mail: dmipak@silk.glas.apc.org}}
\date{23 August 1995}
\maketitle
\begin{center}
       Department of Theoretical Physics, Research Institute\\
       of Applied Physics, Tashkent State  University,\\
       Vuzgorodok, 700095, Tashkent, Republic of Uzbekistan \\
\end{center}

\begin{center}
Classification numbers: 210, 220, 240
\end{center}

\begin{abstract}
     A three-dimensional $q$-Lie algebra of \su  is realized
in terms of first- and second-order differential operators.
Starting from the $q$-Lie algebra one has constructed a left-covariant
 differential calculus on the quantum group. The proposed
 construction is inverse to the standard
Woronowicz approach; the left-invariant vector fields are
introduced as initial objects whereas the differential
 1-forms are defined in a dual manner.
\end{abstract}
 \vskip 50pt

 \newpage
\section{Introduction}
\indent

       Researches of non-commutative geometry of the quantum groups
[1-5]
 led to a series of papers concerned with
 non-equivalent differential calculi on the quantum groups and quantum
spaces. In particular, bicovariant $4D_{\pm}$ and left-covariant
 $3D$ differential calculi on the quantum group \su had been considered
[3-6] .
 For a two-parameter deformed linear
quantum group $GL_{p, q}(2)$ the corresponding covariant differential
calculus on the quantum group and quantum plane was proposed as well
[7 - 9].
 As a rule, the standard Woronowicz approach
\cite{Wor5} is successfully
 used to construct differential calculi on quantum groups.
 In the framework
 of this approach the differential 1-forms are treated as initial
basic elements whereas the vector fields are defined as dual ones.
As a consequence the geometrical content of the vector fields and their
connection with  usual derivatives on the quantum group remain hidden.

       In this paper an explicit representation of the $q$-Lie
algebra of  left-invariant vector fields on the quantum group
\su is proposed in terms of first- and second-order differential
 operators. Using the representation one has constructed the
covariant $3D$ differential calculus on the quantum group.
In some sense, the construction is inverse to the standard
Woronowicz approach. In Section 2 a left-covariant differential calculus
 on  the quantum group \uqtwo is presented. The main commutation
relations are differed from ones considered in \cite{WZ1}.
Here we introduce covariant tensor notations which are very
convenient in constructing various covariant geometrical
objects with a definite \uone charge. Section 3 deals with the
left-invariant vector fields on the quantum group \su.
Defining relations for an exterior differential algebra
complete the differential calculus. In Conclusion
we determine the connection between the \su left-covariant
$3D$ $q$-Lie algebra and the Drinfeld-Jimbo quantum enveloping
 algebra. Main formulas and conventions used throughout the paper
are contained in Appendix.

\section{ Differential calculus on the quantum  group \uqtwo}
\setcounter{equation}{0}
\indent

     Let $x^i, y_i \  (i=1,2)$ be generators (coordinates)
of the function algebra $C^{2}_{q,q^{-1}}$ on the quantum
hermitean vector space $U^2_q$ with an involution
 ${\ast}:{ \stackrel{\ast}{x^{i}} = y_i}$. It is convenient to
parametrize the matrix ${T^i}_j \in \uqtwo $ by the coordinates
$x, y$
\bea
    {T^i}_j=\left(\begin{array}{cc}
                  y^1&x^1                              \label {Tparam}    \\
                  y^2&x^2
                 \end{array} \right)
                             =(y^ix^i).
\eea
 We use the
$R$-matrix formulation of quantum groups
following  Faddeev,
\newline
 Reshetikhin and Takhtajan \cite{Fad1}.
 The main commutation relation for the quantum group generators  $T^i_j$
 has a standard form
\be
   R_{12}T_1 T_2 = T_2 T_1 R_{12}.
\ee

     The parametrization (\ref{Tparam}) was introduced in the harmonic
formalism \cite{Gal1} applied to extended supersymmetric theories and
supergravities. If one imposes the unimodularity constraint
\be
  {\cal D} \equiv {\det}_{q} {T^i}_j = x_i y^i = 1,
\ee
then the variables $(x,y)$ will be just the quantum harmonic functions
$(u^{\pm})$ on the coset  $S^{2}_q  \, \sim \,  {\suu } / {\uone}$ with
 corresponding
 \uone charges $\pm1$:
 \bea
     x^i \equiv u^{+i}, \;\;\;\;\;\;\; y^i\equiv u^{-i}.
 \eea
 For convenience we shall not use the notations $u^{\pm}$
throughout the paper keeping in mind that all geometrical objects
(like coordinates, derivatives, differentials etc.) have definite
\uone charges. All further consideration respects the global
covariance under the action of the group \uone.

   Let us define the derivatives
$\partial_i\equiv \displaystyle{\frac{\partial}{\partial
 x^i}},\;\; \bar\partial^i\equiv \displaystyle
{\frac{\partial}{\partial y_i}}$ on the
quantum group by the formulas
\be
   \partial_ix^j=\gamma_i^j,  \:\:\:\:\:\:\:
  \bar\partial^iy_j=\delta^i_j,
\ee
where $
   \gamma_i^j,\,\,\delta^i_j$ are quantum analogues to classical
 Kronecker symbol. The commutation relations between the
coordinates and derivatives are uniquely defined (up to the
symmetry connected with the exchange $q \rightarrow \qin$):
 \bea
\ba{lr}
R_{12}(\pro_T)_1 (\pro_T)_2 = (\pro_T)_2 (\pro_T)_1 R_{21} ,&
    (\pro_T)^i_j \equiv \left(\begin{array}{cc}
                  \bar \pro_1 &\bar \pro_2\\
                  \pro_1 &\pro_2
                 \end{array} \right)     ,
\ea
\eea
\bea
\ba{lr}
    \pro_ix^k=\gamma_i^k+qY^{nk}_{mi}x^m\pro_n,  &
                        {\bar\pro}^i\bar y_j=\delta^i_j+qy_m\bar\pro^n\hat
                        R^{mi}_{nj},  \\             \label{pro}
    \pro_i y_j=q\hat R^{-1 lk}_{ji} y_k\pro_l,   &
                        \bar \pro^i x^j = \qin \hat R^{ij}_{kl}
                        x^k\bar\pro^l.
\ea
\eea
All these relations are consistent with the quantum group structure.
The definitions (\ref{pro}) do not differ from ones
considered in ref. \cite {WZ1} on principal.
 Our choice is conditioned by a requirement
of manifest \uqtwo
and \uone  covariance. We consider all objects  with upper (lower) indices
 to be  transformed
under the quantum group co-action $ \Delta$ like
classical co-(contra-)variant tensors.
 For instance, a second rank tensor $N_i^j$  will be transformed
as follows
\be
      (N_i^j)^{\prime} = {(T^{\dag})^k}_i {T^j}_l N_k^l .
\ee
   Hereafter the signs $ \bigotimes$ of tensor product are omitted
for simplification in writing formulas. It should be noted
that an alternative way of introducing manifest
  tensor notations is adoption
 of the left $\stackrel{\rightarrow}{\bar\partial^i}$ and right
$\stackrel{\leftarrow}{\partial_i}$ derivatives [10, 11] .

     Let us construct the left-invariant differential operators
$\pro_a \  (a=1,2,3,4)$ on the $U_q(2)$ space:
\bea
\ba{lc}
    \pro_1=x_k \pro^k,   &              \pro_2= y_k\bar\pro^k   ,\\
    \pro_3=x_k\bar\pro^k,  &            \pro_4=y_k\pro^k
\ea
\eea
  The differential operators  $\pro_a$ form a
generalized  $q$-Lie algebra of left-invariant vector fields.
 It is easy to check that the operators $\pro_3,\pro_4$
 maintain the unimodularity condition
 $s \equiv {\cal  D}  -1 =0 $, i.e.
\be
   \pro_{3,4} (sf(x,y))=0.           \label {det}
\ee
Here, $f(x,y)$ is an arbitrary function.  At the same time it
 is impossible to construct a third linearly independent left-invariant
first-order differential operator obeying the condition (\ref {det}).
 Consequently, one cannot realize a  $q$-Lie algebra
on the  \su in terms of first-order differential operators.

    Now we define an exterior differential $\hat d $ by following
relations
\bea
\ba{lc}
    \hat d \equiv d+\bar d,   &       \hat d^2=d^2=\bar d^2=0 \\
    d=dx^i \gamma^{-1 k}_i \pro_k, &   \bar d = dy_i \bar\pro^i
                                                         \label{extd} \\
    d(fg)=df \cdot g+f \cdot dg,   &  \bar d (fg) = \bar d f\cdot g+
                                     f\cdot \bar dg  \\
    d\bar d+\bar d d  =0.   &         \
\ea
\eea
Here a standard Leibnitz rule for the exterior differential
is used. Commutation relations for the basic differential 1-forms,
coordinates and derivatives
can be easily derived in a similar way as in ref.
 \cite {WZ1} taking into account full consistence with
 the complex structure provided by the
 $\ast$-involution. The final relations are collected in Appendix.
The left-invariant Cartan 1-forms $\omega^a (a=1,2,3,4)$
 generate a basis in the space of differential forms
\bea
\ba{lr}
    \omega^1=\qin y_idx^i,     &    \omega^2 = q^2\, x_idy^i,  \\
    \omega^3 =- y_idy^i ,   &    \omega^4 = q\, x_idx^i.
\ea
\eea
The exterior differential defined by relations (\ref{extd})
 can be rewritten as follows
\be
   \hat d =  \frac{1}{\cal D} \omega^a \pro_a.
\ee
Commutation relations for the basic differential
1-forms  $\omega^a$ and the Cartan-Maurer equation are derived
 straightforwardly
using definitions (\ref{extd}) (see Appendix).
 Analysis of the differential
calculus on \uqtwo implies that one cannot
directly reduce it to the differential calculus on the quantum group \su.

\section{Left-invariant vector fields of \su and  \newline
         differential calculus }
\setcounter{equation}{0}
\indent

     To define the left-invariant vector fields on the quantum
group \su \ we shall use the $U_q(2)$-covariant differential
calculus. Let us introduce the next notations for the
 left-invariant 1-order differential operators in correspondence
with notations of the classical harmonic approach \cite {Gal1}
\bea
    D^{++}\equiv x_i\bar\pro^i,   &   D^{--}\equiv -y_i\pro^i.
\eea
The action of the operators $D^{++},\, D^{--}\,$ on the coordinates
$(x, y)$
have simple properties:
\bea
\ba{lr}
    D^{++}x^i=0,  &    D^{--}x^i=y^i,       \label {def1}  \\
    D^{++}y_i=x_i,  &  D^{--}y_i=0.
\ea
\eea
 A Leibnitz rule for these operators is simplified when
acting on the functions with  definite
$U(1)$ charges
\be
   D^{\pm\pm}(f^{(m)} g^{(n)}) = (D^{\pm\pm}f^{(m)}) g^{(n)}+
                         q^{-m}f^{(m)}( D^{\pm\pm}g^{(n)}).
                                             \label {def2}
\ee

     As it is mentioned above, one cannot construct
a $q$-analogue for the classical \uone \  generator $D^0$ in terms of
first-order differential operators.
Nevertheless, the interesting feature of non-commutative geometry
is that the quantum \uone generator $D^0$
does exist. The operator is realized as a left-invariant
 second-order differential operator
\be
   D^0\equiv -x_i\pro^i-q^2 y_i\bar\pro^i +
                    (1-q^2)x_iy_k\bar\pro^k\pro^i.
\ee
It is not hard to check that the
operator $D^0$ has eigenfunctions $f^{(n)}$ with eigenvalues
corresponding to \newline $q$-generalized \uone charge $(n)$:
\bea
   D^0f^{(n)} ={\{n\}}_q f^{(n)},             \label{def3}            \\
{\{n\}}_q \equiv    \dfrac {1-q^{-2n}}
                          {1-q^{-2}}.     \nn
\eea
Due to that property (\ref{def3}) the algebra of functions $f^{(n)}$
with a definite \uone charge is just the function algebra on
a quantum sphere $S^2_q = \su /U(1)$. The Leibnitz rule for the
operator $D^0$ has the next form
\be
   D^0(f^{(m)}g^{(n)}) = (D^0f^{(m)})g^{(n)} + q^{-2m}f^{(m)}
                                     D^0g^{(n)} .          \label {def4}
\ee

By direct checking one can verify that the operators
$D^{\pm\pm 0}$ satisfy a generalized $q$-Lie
algebra of \su \ \cite {Schi1}
\bea
\ba{lr}
  {[D^0,D^{++}]}_{\dpst {q^{-4}}} = {\{ 2\}} _q D^{++},&   \\
   {[D^0, D^{--}]}_{\dpst {q^4}} = {\{ -2\}} _q D^{--}, &     \label{alg}\\
   {[D^{++}, D^{--}]}_{\dpst {q^2}} = D^0,             &
\ea
\eea
here, $[A,B]_{\dpst {q^s}} \equiv AB-q^s BA$.

      It should be noted that the algebra (\ref {alg})  is valid
irrespective of whether one imposes the constraint ${\cal D}=1$.
An important property of the operators $D^{\pm\pm 0}$
is conservation of the unimodularity constraint
\be
   D^{\pm\pm0}({\cal D}\: f(x,y))\cong 0.
\ee
The last relation allows to construct the differential calculus on
the \su in a consistent manner. Note, that the braiding matrix
corresponding to the $q$-Lie algebra (\ref {alg}) is unitary
 and the generalized Jacobi
 identity is available
\bea
    [D^0, [D^{++}, D^{--}]_{\dpst {q^2}}]
 + [D^{++},[D^{--},D^0]_{\dpst {q^{-4}}}]_{\dpst {q^{-2}}}  \nn \\
            +q^2[D^{--},[D^0,D{++}]_{\dpst {q^{-4}}}]_{\dpst {q^{-2}}}
             \equiv 0.
\eea
It is easy to check another relation which is similar to Jacobi
identity
\bea
    [D^0, [D^{++}, D^{--}]] + [D^{++},[D^{--},D^0]_{\dpst {q^{-4}}}]
_{\dpst {q^6}}                 \nn   \\
            + q^6[D^{--},[D^0,D{++}]_{\dpst {q^{-4}}}]_
{\dpst {q^{-6}}} \equiv 0,
\eea
Using properties of the differential operators $D^{\pm \pm0}$
we can define a covariant algebra of left-invariant vector fields
$\nabla^{\pm\pm0}$ on the quantum group \su \ in axiomatic way and then
introduce the basic differential Cartan
 1-forms as dual objects. Let us define
the left-invariant vector fields $\nabla^{\pm\pm0}$ by the same
relations (\ref {def1}, \ref {def2}, \ref {def3}, \ref {def4})
 that the operators $D^{\pm\pm0}$ obey
with  only exchanging $D^{\pm\pm0} \rightarrow \nabla^{\pm\pm0}$.
The basic left-invariant differential 1-forms $\omega^{\pm\pm0}$
are then defined as dual objects to vector fields
$\nabla^{\pm\pm0}$
\bea
   \omega^{++}(\nabla^{--}) =1, \;\;\;  \omega^{--}(\nabla^{++})=q,\;\;\;
                              \omega^0(\nabla^0)=1.
\eea
An exterior differential on \su \ is defined in a standard manner
with a usual Leibnitz rule
\bea
\ba{lr}
   \delta \equiv \omega^{++} \nabla^{--} +
    \omega^{--}\nabla^{++} + \omega^0\nabla^0,  &      \\
   \delta (fg) = \delta f\cdot g+ f\cdot\delta g, &      \\
   \delta^2 =0,                                &
\ea
\eea
where the $f, g $ -- are arbitrary functions on the quantum
 group  \su. Using these formulas it is not difficult to obtain
all commutation relations for the differential 1-forms
 $\omega^{\pm\pm0}$ and corresponding
Cartan-Maurer equations
\bea
\ba{lr}
   {(\omega^{\alpha})}^2 =0, \;\;\;  \alpha =(++,--,0), &
                         \op\om=-{\displaystyle \frac {1}{q^2}}
                                              \om\op,    \\
    \omega^{\pm\pm} f^{(m)} = q^m f^{(m)}\omega^{\pm\pm},  &
                                       \op\oo =-{\displaystyle
                          \frac{1}{q^4}}\oo\op,       \\
   \oo f^{(m)} = q^{2m}f^m\oo,    &    \om\oo=-q^4 \oo\om,  \\
  \delta \op= {\{-2\}}_q \op\oo,    &                   \\
  \delta \om = {\{2\}}_q \om\oo,    &                   \\
  \delta \oo = \op\om.    &
\ea
\eea
All these relations are consistent with the unimodularity
 condition ${\cal D}=1$.
 An exterior algebra of the differential forms is defined
straightforwardly. The final construction of the
 covariant differential calculus on the quantum group \su \ presented here
agrees with  one considered in  Woronowicz approach
[3, 8] .

\section{Conclusion}
\setcounter{equation}{0}
\indent

      Now we shall determine the explicit connection between the
generalized $q$-Lie algebra (\ref {alg}) and the quantum
enveloping Drinfeld-Jimbo algebra ${\got U}_q (su(2))$. For
this purpose one considers the differential operators $\mu ,\ \nu$
\cite {Ogi1} on the quantum group $U_q(2)$:
\be
   \mu = 1+(q^2-1) y_i{\bar {\pro}}^i, \;\;\;\;\;
   \nu = 1+ (1- \frac{1}{q^2}) x_i \pro^i.
\ee
These operators have simple commutation relations with the
operators $D^{\pm\pm 0}$. For instance, we have the following
formulae
containing the
 operator
 $ \mu $
\bea
\ba{lr}
   \mu D^{--} = q^2 D^{--} \mu, \;\;\;\; & \mu D^{++} =
\dfrac {1}{q^2} D^{++}\mu,\\
   \mu D^0 = D^0 \mu,            &   \mu\nu=\nu\mu  .
\ea
\eea
Equations for the $ \nu$ operator have a similar form.

     Let us define new operators $ \dpp ,\dmm , \docal $
 multiplying the $D^{\pm\pm 0}$
by corresponding factors
\bea
\ba{lr}
    \dpp = \mu^{-\frac {1}{2}} D^{++},& \dmm = \nu^{-\frac{1}{2}} D^{--} ,\\
    \docal = \qin \mu \nu D^0  \equiv {[\pro^0]}_q.    &
\ea
\eea
It is easy to verify that the operators ${\cal D}^{\pm\pm 0}$ form just
the Drinfeld-Jimbo quantum enveloping algebra ${\got U}_q (su(2))$
\bea
 &  [\pro^0, \dpp ] = 2\dpp ,    &   [\pro^0, \dmm ] = -2\dmm , \nonumber  \\
 &  [\dpp ,\dmm ] = {[\pro^0]}_q.  &
\eea
The operator $\docal$ counts the $q$-generalized $U(1)$ charge when
 acting on the functions with a definite $U(1)$ charge $(m)$
\be
   \docal f^{(m)} = [m]_q f^{(m)}.
\ee
Observe that the way of constructing the Drinfeld-Jimbo algebra
from the $q$-Lie algebra (\ref {alg}) is by no means unique.

        We confine our consideration to a simple case of the
quantum group \su. However, one should expect that similar
representation of the $q$-Lie algebra in terms of
differential operators hold for other quantum groups. The
$q$-Lie algebra (\ref {alg}) has a natural geometrical
origin in our approach. It turns out to be closely connected
with a gauge covariant differential algebra of \su in
constructing the non-standard Leibnitz rule. These questions
will be considered in a separate paper elsewhere.

\bigskip
{\bf Acknowledgments}
$$~$$
    Author would like to thank B. Zupnik,
 A. Isaev and Ch. Devchand for
useful discussions and interest to work.
$$~$$

\appendix
\section{Appendix}
\setcounter{equation}{0}
\indent

       We use the next notations for the invariant \su \ tensors
\bea
    \vareps^{ij} = \left ( \begin{array}{cc} 0  &  1  \\
                                         -q &  0   \end{array} \right ),
                      \;\;\;\;\;\;\; &
                         \gamma_{i}^{j} = \left ( \begin{array}{cc}
                        q  &  0  \\
                        0 &\dfrac {1}{q} \end{array} \right )  \nonumber   \\
    \vareps_{ij} = \left ( \begin{array}{cc}  0  &  1  \\
                                   -\dfrac {1}{q}&0   \end{array} \right ),
                       \;\;\;\;\;\;\; &
                        {\gamma ^{-1}}^j_i= \left ( \begin{array}{cc}
                        \dfrac {1}{q}  &  0  \\
                        0        &  q    \end{array} \right )    \\
     \vareps_{ik} \vareps^{jk} = \vareps_{ki} \vareps^{kj} = \delta^j_i ,
                       \;\;\;\;\;\;\;    &
                  \vareps_{ik} \vareps^{kj} = -{\gamma}_i^j  \nonumber \\
     \vareps_{ki} \vareps^{jk} = -{\gamma ^{-1}}^j_i.    &  \nonumber
\eea
The invariant metric $\vareps_{ij}$ is used to raise  and lower the
\su \ indices  as follows
\bea
    A_i = \vareps_{ij} A^j,   &    A^i = A_j \vareps^{ji}
\eea

     The $R$-matrix and auxiliary matrices $X, Y$ are defined as
in \cite {WZ1}
\bea
   {\hat R}^{ji}_{kl} =
 R^{ij}_{kl} = \delta ^i_k\delta^j_l(1+(q-1)\delta^{ij})
      +(q-{\dfrac {1}{q}})\delta^j_k\delta^i_l\theta (i-j),      &   \\
\ba{ll}
    X^{ri}_{sj} = {\hat R}^{ir}_{js} q^{2(r-j)} =
                  {\hat R}^{ir}_{js} q^{2(s-i)},               &
                     X^{ri}_{sj}{\rmm}^{jk}_{il} = \delta^r_l
                                             \delta^k_s, \\
    Y^{ri}_{sj} = {\rmm}^{ir}_{js} q^{2(s-i)} =
                 {\rmm}^{ir}_{js} q^{2(r-j)},                &
               Y^{jk}_{il} {\hat R}^{ln}_{km} = \delta^n_i\delta^j_m .
\ea                                                              \nn
\eea

    The main commutation relations in a case of the quantum group
$U_q(2)$ have the following form
\bea
\ba{ll}
     R_{12}dT_1dT_2=-dT_2dT_1R^{-1}_{12} ,  &  \\
dx_1 x_2 = q R_{21} x_2dx_1,  &
                   \pro_j dx^i = q X^{ki}_{lj}dx^l\pro_k,        \\
dy_1 y_2 = \qin R^{-1}_{12} y_2dy_1, &
             \bar{\pro}^idx^j = \qin X^{ji}_{lk} dx^k\bar{\pro}^l, \\
dy_1 x_2 =  R_{21} x_2dy_1,  &
             {\bar \pro}^idy^j = \qin X^{ji}_{lk} dy^k {\bar \pro}^l,\\
dx_1 y_2 =  R^{-1}_{12} y_2dx_1, &
    {\pro}_i{dy}_j = q{\rmm}^{lk}_{ji} {dy}_k {\pro}_l.      \\

\ea
\eea

     The Cartan-Maurer equations for the quantum group \uqtwo
can be written as follows
\bea
\ba{lr}
   d\omega_1 = {\dpst \frac {1}{\cal D}}(\omega_3\omega_4 +
                                             \omega_2 \omega_1), &
   d\omega_2 =-{\dpst \frac {1}{\dpst \cal D}} (q^2\omega_4\omega_3 +
                 \omega_1\omega_2),   \\
   d\omega_3 = -{\dpst \frac {1+q^2}{q^4 \cal D}} \omega_3\omega_2,     &
   d\omega_4 = {\dpst \frac {1+q^2}{\cal D}} \omega_4\omega_1.
\ea
\eea

\end{document}